\documentclass[12pt]{article}
\usepackage[cp1251]{inputenc}
\usepackage[dvips]{graphics,graphicx}
\usepackage[english]{babel}
\usepackage{amsfonts,amssymb, amsmath, theorem}

\def\RR{\hbox{I\kern-.2em\hbox{R}}}

\newtheorem{Def}{{\mbox{$\;\;\;\;\;\,$}}Definition}[section]
\newtheorem{Th}{{\mbox{$\;\;\;\;\;\,$}}Theorem}[section]

\newtheorem{uess}{{\mbox{$\;\;\;\;\;\,$}}Lemma}[section]

\newtheorem{Rem}{{\mbox{$\;\;\;\;\;\,$}}Remark}[section]

\usepackage{graphicx}
\begin{document}
\author{P. Amster
   \\Universidad de Buenos Aires and CONICET\\
      Departamento de Matem\'atica - FCEyN \\
Ciudad Universitaria, Pab. I\\
       1428 - Buenos Aires, Argentina \\
   email:  {\tt{pamster@dm.uba.ar}}\\\\
    L. Berezansky
   \\
      Department of Mathematics\\
   Ben-Gurion University of Negev  \\
   Beer-Sheva 84105, Israel  \\
   email:  {\tt{brznsky@cs.bgu.ac.il}} \\\\
  L. Idels
         \\
   Department of Mathematics \\
   Vancouver Island University  \\
   900 Fifth St. Nanaimo, BC, Canada \, V9S5S5   \\
   email:  {\tt{lev.idels@viu.ca}}
}
\newcounter{tally}
\date{}

\title{Stability of Hahnfeldt Angiogenesis Models with Time Lags}
\maketitle \vfill{}

\begin{abstract}
Mathematical models of angiogenesis, pioneered by P. Hahnfeldt, are under study.
To enrich the dynamics of three models, we introduced biologically motivated time-varying delays. All models under study belong to a special class of nonlinear nonautonomous systems with delays. Explicit conditions for the existence of positive global solutions and the equilibria solutions were obtained. Based on a notion of an M-matrix, new results are presented for the global stability of the system and were used to prove local stability of one model. For a local stability of a second model, the recent result for a Lienard-type second-order differential equation with delays was used. It was shown that models with delays produce a complex and nontrivial dynamics. Some open problems are presented for further studies.
\end{abstract}
{\bf. Keywords} Angiogenesis, Nonlinear nonautonomous delay differential equations, Global and local stability, Equilibria, M-matrix, Lienard equations. MSC 2000: 34K06, 34K20, 34K60.
\section{Introduction}
Angiogenesis, the generation of new blood vessels, is thought to be necessary for tumor growth and metastasis  \cite{Folk2}--\cite{Folk3}. This process belongs to a general family of tumor-immune interactions. Some recent studies, including mathematical models of angiogenesis dynamics, are presented in \cite{A},  \cite{Pil}, \cite{O2}, \cite{O4}-- \cite{Folk3}, \cite{Hah1}--\cite{Hah4}, \cite{KK},  \cite{Sa} and \cite{Sw}.\\
To incorporate the spatial effects of the diffusion factors that stimulate and inhibit angiogenesis, the following two-compartmental model  for cancer cells and vascular endothelial cells was developed by P. Hahnfeldt in  \cite{Hah1} (see, also \cite{Sa}).

\begin{equation}\label{HFolk}
\begin{array}{ll}
\displaystyle\frac{dx}{dt}=\alpha x(t)\ln \frac{K(t)}{x(t)} \\[3mm]
\displaystyle\frac{dK}{dt}=-\mu K(t)+ S(x(t),K(t))-I(x(t), K(t))- c(t) K(t),
\end{array}
\end{equation}
 where $x(t)$ is the tumor mass and $K(t)$ is a variable carrying capacity, that is defined as the effective vascular support provided to the tumor as reflected by the size of the tumor potentially sustainable by it.
According to P. Hahnfeldt, a stimulator/inhibitor tumor growth dynamics, described by system \eqref{HFolk}, should provide a time dependent carrying capacity under angiogenic control and include the distinct mechanisms for angiogenic stimulation and inhibition.
The dynamics of the second equation is a balance between stimulatory and inhibitory effects:
first term is the loss of functional vasculature; the second term corresponds the stimulatory capacity of the tumor.
The third term reflects endogenous inhibition  by either neutralizing endothelial cell growth factors or inhibition endothelial cell proliferation.
These inhibitors are released through the tumor surface (scaling the tumor volume to its surface area); thus the major assumption for the Hahnfeldt models is
\begin{equation}
I(x,K)/S(x,K)\sim K^{\theta}x^{\mu} \nonumber,
\end{equation}
where $\theta+\mu = 2/3$.\\
 The last term in model \eqref{HFolk} represents inhibition of tumor vasculature due to administered therapy. The following pharmokinetic differential equation
$$\frac{dc}{dt}=v(t)-qc(t)$$ might be used to model chemo- or radio-therapy mechanisms.
Here $v(t)$ is a dose given and $q$ is a per capita decay rate of the drug once it is injected.
If we assume that the drug kills all types of cells \cite{Hah4}, then
\begin{equation}\label{Ang3}
\begin{array}{ll}
\displaystyle\frac{dx}{dt}=\alpha x(t)\ln \frac{K(t)}{x(t)}-p(t) x(t) \\[3mm]
\displaystyle\frac{dK}{dt}=-\mu K(t)+ S(x(t),K(t))-I(x(t), K(t))- c(t) K(t).
\end{array}
\end{equation}
To model processes in nature it is frequently required to know system states from the past, i.e., models incorporating memory. Depending on the phenomena under study the after-effects represent duration of some hidden processes, for example, time lags of transit through one state to another; transit time through compartments, or time lags associated with the growth rates (cell division/differentiation time).
In any cell growth some cells are inactive and, once activated, the cell division is not instantaneous.
There are delays in cell division. For example, cells can be engaged in an active cell cycle and might divide with a fairly regular periodicity. On the other hand, cells can be in a resting state (i.e., having dropped out of the cell cycle because of lack of positive growth signals) and will experience a delay in resuming the cell cycle at the same place that they left. In human cells this delay  is usually about 8 hours. One might imagine that in addition to carrying oxygen and nutrients to tumor cells in a tumor mass, vascularization may also deliver growth hormones (signals) that would stimulate otherwise quiescent cells to reenter the cell division cycle and thereby contribute to an increase in tumor mass.
The Warburg effect on angiogenesis and metastasis is a different mechanism with the after-effects (see, for example, \cite{W}). Essentially cancer cells continue to use glycolysis much more than regular cells, even when there is plenty of oxygen around. This produces an acidic environment (by the production of lactic acid from glycolysis) which in turn stimulates angiogenesis. Presumably there is a lag time required to bring the tumor environment to the optimal pH for endothelial cell growth.
In terms of inhibition, there could be time lags between cells receiving an inhibitory signal and actually dropping out of the cell cycle. Usually cells will continue through their cycle after receiving a signal and all stop at the same spot on the cycle (called a restriction point).
The inclusion of explicit time lags in the model allows direct reference to experimentally measurable and/or controllable cell growth characteristics (e.g., time required to perform the necessary divisions). In general, models with delays produce a complex and nontrivial dynamics: it switches stable trajectories into unstable cycles or periodic oscillations. In cancer therapy, stability switching is a very important issue in the design of a drug protocol (see, for example, \cite{And}, \cite{Bab}, \cite{O3} and \cite{Xu}).  \\

If we assume that the tumor cells enter the mechanisms of angiogenic stimulation and inhibition with some delays $h(t)\leq t$, then model \eqref{HFolk} has two alternative forms:\\
\textit{Model 1.}
\begin{equation}\label{Ang0}
\begin{array}{ll}
\displaystyle\frac{dx}{dt}=\alpha x(t)\ln \frac{K(t)}{x(t)}-p(t)x(t) \\[3mm]
\displaystyle\frac{dK}{dt}=\beta x(h(t))-\gamma K(t) -\delta x^{2/3}(h(t)) K(t)- c(t) K(t).
\end{array}
\end{equation}
\textit{Model 2.}
\begin{equation}\label{Ang00}
\begin{array}{ll}
\displaystyle\frac{dx}{dt}=\alpha x(t)\ln \frac{K(t)}{x(t)}-p(t)x(t) \\[3mm]
\displaystyle\frac{dK}{dt}=\beta K(t)-\gamma K(t)-\delta x^{2/3}(h(t)) K(t)- c(t) K(t).
\end{array}
\end{equation}
Note that a logistic-type model with Richards nonlinearity\\
\textit{Model 3.}
\begin{equation}\label{AngLog}
\begin{array}{ll}
\displaystyle\frac{dx}{dt}=\alpha x(t)\left(1-\left[\frac{x(t)}{K(t)}\right]^{m}\right )-p(t)x(t)  \\[3mm]
\displaystyle\frac{dK}{dt}= \beta x(h(t))-\gamma K(t)-\delta x^{2/3}(h(t)) K(t)- c(t) K(t).
\end{array}
\end{equation}
could also be used for modelling tumor growth dynamics \cite{Hah1}. Here $m>0$ and $m \neq 1$ is a constant
that drops an unnatural symmetry of the classical logistic curve ($m=1$).\\
Note that all models without time lags  were studied in \cite{O2}, \cite{O4}, \cite{Hah1} and \cite{Sa}.\\

Systems \eqref{Ang0}--\eqref{AngLog} belong to a wide class of the delay differential equations
\begin{equation}
\label{12A}
\frac{dx}{dt}= A(t)x(t)+F(t, x(t),x(h(t))), ~t\geq t_0,
\end{equation}
where $x \in {\mathbb R}^{n}$, $A \in {\mathbb R}^{n\times
n}$ and $F : {\mathbb R}_{+}\times {\mathbb R}^{n}\times {\mathbb R}^{n}\rightarrow {\mathbb R}^{n}$
is a nonlinear and continuous vector function.\\
Stability analysis of the delay differential equation \eqref{12A} is a well-trodden area,
however, some existing results  rely on restrictive conditions, e.g., strict monotonicity and boundedness of the functions and operators involved, continuity of the parameters \cite{Ber3}, \cite{Dib1},  \cite{Gil2}, \cite{G2}, \cite{Lev}, \cite{Kir} and \cite{Mur}. For example, models under study possess non-Lipchitz nonlinearity, thus stability analysis for models \eqref{Ang0}--\eqref{AngLog} requires new tools and approaches. Moreover, since the models in the paper are based on nonautonomous equations, the application of the traditional methods  for autonomous equations, such as Laplace transforms or quasipolynomials, is questionable.\\

The purpose of this paper is two-fold: obtain new results in the qualitative theory of delay differential equations; apply the results to the analysis of the models of angiogenesis dynamics.
Firstly, for some general class of nonlinear nonautonomous systems with delays we obtained conditions for the existence of a global attractor; and based on that result, we proved local stability for Model l. Compared with Model 1, Model 2 has different qualitative features, thus for its analysis, we used the theorems recently obtained for a Lienard-type second-order differential equation with delays \cite{Ber2}. In the Discussion we pose some open problems.

\section{Global Stability Criteria}
We will begin by examining the global stability of a general class of the following nonlinear nonautonomous systems with delays.
\begin{equation}
\label{nonlinear}
\dot{x_i}(t)=-a_{ii}(t)x_i(t)+\sum_{j\neq i}a_{ij}(t)x_j(t)
+\sum_{j=1}^n f_{ij}(t,x(h_{ij}(t))),~i=1,\dots, n,
\end{equation}
where $|f_{i,j}(t,u)|\leq b_{ij}(t)|u|$, $a_{ij}(t) ~ and ~ b_{ij}(t)$ are measurable essentially bounded functions, $f_{ij}(t,\cdot)$ is a continuous function, $f_{ij}(\cdot, u)$ is a measurable locally essentially bounded function and $h_{ij}(t)\leq t$ are measurable functions.
Along  with equations (\ref{nonlinear}) we set the initial condition
for each $t_0\geq 0$
\begin{equation}\label{in}
x_i(t)=\varphi_i(t), ~ t\leq t_0, i=1,\dots,n,
\end{equation}
where $\varphi_i$ is a continuous function.\\
A solution $X(t)=\{x_1(t),\dots, x_n(t)\}^T$ of problem  (\ref{nonlinear})-(\ref{in})  is a vector function,
locally absolutely continuous for $t\geq t_0$, that satisfies almost everywhere equation (\ref{nonlinear}) on this interval and initial condition (\ref{in})
for $t\leq t_0$.  A unique global solution of problem (\ref{nonlinear})-(\ref{in}) exists.

{\bf Definition.}
{\em Matrix $A$ is called an {\em $M$-matrix} if $a_{ij}\leq 0, i\neq j$ and
one of the following equivalent conditions holds:

-there exists a nonnegative inverse matrix $A^{-1}\geq 0$.

-the main minors of matrix $A$ are positive numbers.}\\

Let $A_{ij}=\sup_{t\geq t_0} |a_{ij}(t)| ~ and~ B_{ij}=\sup_{t\geq t_0} |b_{ij}(t)|$.

\begin{Th}\label{Sta}
Suppose there exist $a_i>0, \tau>0$  and $ t_0\geq 0$ such that $\inf_{t\geq t_0}a_{ii}(t)\geq a_i$ and  $t-h_{ij}(t)\leq \tau$. The  matrix $B=\{b_{ij}\}$ with entries $b_{ii}=a_{i}-B_{ii}$  and $  b_{ij}=-A_{ij}-B_{ij}$ for $  i\neq j$
is an M-matrix. Then for any solution  $X(t)=\{x_1(t),\dots, x_n(t)\}^T$ of equation (\ref{nonlinear})
 $\lim_{t\rightarrow\infty} X(t)=0$.
\end{Th}

{\bf Proof.}
After substitution
$x_i(t)=e^{-\lambda (t-t_0)}y_i(t), ~t\geq t_0$, where $0<\lambda < \min_{i} a_{i} $, system (\ref{nonlinear}) has a form
$$
\dot{y}_i(t)=-(a_{ii}(t)-\lambda)y_i(t)+\sum_{j\neq i}a_{ij}(t)y_j(t)+\sum_{j=1}^n e^{\lambda (t-t_0)}f_{ij}(t,e^{-\lambda (h_{ij}(t)-t_0)}y_j(h_{ij}(t))).
$$
Hence
$$
y_i(t)=e^{-\int_{t_0}^t [a_{ii}(s)-\lambda]ds}x_i(t_0)+
\int_{t_0}^t e^{-\int_{s}^t [a_{ii}(\zeta)-\lambda]d\zeta }\left[\sum_{j\neq i}a_{ij}(s)y_j(s)\right.
$$$$
+\left.\sum_{j=1}^n e^{\lambda (s-t_0)}f_{ij}(s,e^{-\lambda (h_{ij}(s)-t_0)}y_j(h_{ij}(s)))\right]ds.
$$
Then
$$
|y_i(t)|\leq e^{-\int_{t_0}^t [a_{ii}(s)-\lambda]ds}|x_i(t_0)|+
\int_{t_0}^t e^{-\int_{s}^t [a_{ii}(\zeta)-\lambda]d\zeta }\left[\sum_{j\neq i}|a_{ij}(s)||y_j(s)|\right.
$$$$
+\left.\sum_{j=1}^n |b_{ij}(s)|e^{\lambda (s-h_{ij}(s))}|y_j(h_{ij}(s))|\right]ds
$$
$$
\leq |x_i(t_0)|+ \int_{t_0}^t e^{-\int_{s}^t [a_{ii}(\zeta)-\lambda]d\zeta }(a_{ii}(s)-\lambda)\left[\sum_{j\neq i}\frac{|a_{ij}(s)|}{a_{ii}(s)-\lambda}|y_j(s)|\right.
$$
$$
+\left.\sum_{j=1}^n \frac{e^{\lambda \tau}|b_{ij}(s)|}{a_{ii}(s)-\lambda}|y(h_{ij}(s))|\right]ds.
$$
 Set $y_i^b=\max_{t_0-\tau\leq t\leq b}|y_i(t)|~ and~ Y^b=\{y_1^b,\dots,y_n^b\}^T$. Hence
$$
y_i^b\leq |x_i(t_0)|+ \sum_{j\neq i} \frac{A_{ij}}{a_{i}-\lambda}y_j^b+\sum_{j=1}^n \frac{e^{\lambda \tau}B_{ij}}{a_{i}-\lambda}y_j^b
$$
for $t_0\leq t\leq b$.
We define the matrix $C(\lambda)=\{c_{ij}(\lambda)\}$ with entries
$$
c_{ii}(\lambda)=1-\frac{e^{\lambda \tau}B_{ii}}{a_{i}-\lambda}~ and ~ c_{ij}(\lambda)=-\frac{A_{ij}+e^{\lambda \tau}B_{ij}}{a_{i}-\lambda},~ i\neq j.
$$
Hence  $C(\lambda) Y^b\leq |X(t_0)|$ for $t_0\leq t\leq b$.
We have $\lim_{\lambda\rightarrow 0}C(\lambda)=C(0)$. Since matrix $B$ is an M-matrix, then $C(0)$  is  also an M-matrix; thus the main minors of  $C(0)$ are positive. If for some values of parameters the elements of a matrix are continuous functions then the determinant
of this matrix is a continuous function. For some $\lambda >0$ all main  minors of $C(\lambda)$ are positive, the latter implies  that this matrix is an M-matrix.
Assume that parameter is fixed  $\lambda=\lambda_0$.
 Thus for $Y^b$ there is an \textit{a priory} estimate $Y^b\leq C^{-1}(\lambda_0)|X(t_0)|$ where the right-hand side does not depend on $b$; therefore $Y(t)$ is a bounded function. Finally, $X(t)=e^{-\lambda(t-t_0)}Y(t)$, then $\lim_{t\rightarrow\infty}X(t)=0$.
The theorem is proven.

\section{Existence of Global Solutions}
In what follows, we assume
  $0\leq t-h(t)\leq \tau$,
  $x(t)=\phi (t)~ for ~ t<t_{0}$  $x(t_{0})=x_{0}$ and $K(t_{0})=K_0$,
  where $p (t)~ and ~ c(t)$ are measurable essentially bounded functions, $h(t)$ is a measurable function and $\phi$ is a continuous function. \\
 \begin{Def}  Any solution of problems (\ref{Ang0})-- (\ref{AngLog}) is a locally absolutely continuous function if it satisfies the equation almost everywhere for $t>t_0$
  and the initial conditions for $t\leq t_0$.\end{Def}

 \begin{Th}\label{1}
  Suppose that $t-h(t)\geq \tau_0>0, ~ \phi (t)\geq 0, ~ x(t_0)>0 ~and~ K(t_0)>0$. Then systems (\ref{Ang0})-- (\ref{AngLog}) have unique global solutions $(x(t), K(t))$
 positive for $t\geq t_0$.
  \end{Th}
  {\bf Proof. }
 Without loss of generality, we prove the theorem  for system (\ref{Ang0}). The proof for models 2 and 3 is similar. Suppose for simplicity that $t_0=0$.
 Consider the second equation in  system (\ref{Ang0}) for $t\in [0,\tau_0]$. This linear equation has a form
 $$
 \dot{K}(t)=a(t)K(t)+b(t),~ K(0)>0,~  b(t)\geq 0.
 $$
 Hence $K(t)>0$ for $ t\in [0,\tau_0]$.
Consider the first  equation in  system (\ref{Ang0}) for $t\in [0,\tau_0]$. This is a nonlinear ordinary differential
 equation for the function $x(t)$ with a known positive function $K(t)$.
 Since $x(0)>0$, then there exists a local solution of this equation.
On the interval of the existence of this solution we have
$$
x(t)=x(t_0)e^{\int_{0}^t (\alpha\ln \frac{K(s)}{x(s)}-p(s))ds};
$$
hence on this interval the solution is positive.
Suppose that the maximum interval of the existence of the solution of this equation is $[0,t_0)$ for $t_0<\tau_0$.
Since $x(t)>0$ then $\lim_{t\rightarrow t_0}x(t)=+\infty$.  Therefore there exists $0<t_1<t_0$ such that $$x(t)>\max_{0\leq t\leq \tau_0} K(t)~ for~ t\in [t_1,t_0).$$
We have $\ln \frac{K(s)}{x(s)}<0,  t\in [t_1,t_0)$. Then $$ \dot{x}(t)\leq -p(t)x(t) ~ for ~ t\in [t_1,t_0).$$ Finally, $$0<x(t)\leq x(t_1)e^{\int_{t_1}^{t_0} |p(s)|ds}~ for~ t\in [t_1,t_0).$$
It contradicts the assumption that $\lim_{t\rightarrow t_0}x(t)=+\infty$. Therefore, there exists a positive solution of this equation for $t\in[0,\tau_0]$.
Similarly, we apply the same procedure on the intervals $[\tau_0, 2\tau_0], [2\tau_0, 3\tau_0], \dots$ and  obtain the global positive solution for system  (\ref{Ang0}).
The theorem is proven.
Following the steps in this Theorem, it is straightforward to check the existence of the global solutions of models (\ref{Ang00}) and (\ref{AngLog}).\\
In view of Theorem \ref{1} we assume that there exists $0<\tau_1<\tau_0$ such that $\tau_1\leq t-h(t)\leq \tau_0$, e.g., $h(t)=t-\tau$.\\

\section{ Models with Constant Rate of Infusion}
Standard chemo- and radio- therapies are typically administered in a constant dose scheduling, i.e., $c(t)=c_0$ and $p(t)=p_0$.\\
\subsection{Stability Analysis for Model 1}

\begin{equation}\label{Ang}
\begin{array}{ll}
\displaystyle\frac{dx}{dt}=\alpha x(t)\ln \frac{K(t)}{x(t)}-p_0x(t) \\[3mm]
\displaystyle\frac{dK}{dt}=\beta x(h(t))-\gamma K(t) -\delta x^{2/3}(h(t)) K(t)- c_0 K(t).
\end{array}
\end{equation}
Then system
\begin{equation}
\begin{array}{ll}
\alpha x\ln \frac{K}{x}-p_0 x=0 \\[3mm]
\beta x-\gamma K-\delta x^{2/3} K-c_0 K=0 \nonumber,
\end{array}
\end{equation}
has a unique positive equilibrium point
\begin{equation} \label{Equi}
\begin{array}{ll}
x^{*}=\left(\frac{\eta -\gamma-c_0}{\delta}\right)^{\frac{3}{2}}\\[2mm]
K^{*}=x^{*}e^{\frac{p_{0}}{\alpha}}
\end{array}
\end{equation}
where $\eta=\beta e^{-\frac{p_{0}}{\alpha}}$,
provided that \begin{equation} \label{Cond} \beta >(\gamma +c_0)e^{\frac{p_0}{\alpha}}.\end{equation}
Let introduce new variables $u(t)=\ln\frac{x}{x^{*}}$ and $v(t)=\ln\frac{K}{K^{*}}$, then system \eqref{Ang} has the following exponential form with the trivial equilibrium
\begin{equation}\label{Ang2}
\begin{array}{ll}
\displaystyle\frac{du}{dt}= -\alpha  u(t)+\alpha v(t)\\[3mm]
\displaystyle\frac{dv}{dt}=\eta e^{u(h(t))-v(t)} -\gamma -c_0+(\eta-\gamma-c_0 ) e^{^{\frac{2u(h(t))}{3}}}.
\end{array}
\end{equation}
At $(0 , 0)$ a linearization of system \eqref{Ang2} has the following form
 \begin{equation}
\begin{array}{ll}
\displaystyle\frac{du}{dt}=-\alpha  u(t)+\alpha v (t) \\[3mm]
\displaystyle\frac{dv}{dt}=\frac{\eta+2\gamma+2c_{0}}{3}u(h(t)) -\eta v(t) \nonumber.
\end{array}
\end{equation}

 \begin{Th}\label{Th4.1}
 Let $t-h(t)\leq \tau$ and condition \eqref{Cond} holds. Then positive equilibrium $(x^{*},K^{*})$ of system \eqref{Ang} is locally asymptotically stable.
  \end{Th}
{\bf Proof.} To apply Theorem \ref{Sta} we set a matrix

 $$B=
\left[\begin{array}{ll}
	\alpha&-\alpha\\
	-\frac{\eta+2\gamma+2c_{0}}{3}& \eta
\end{array}\right].$$
Clearly, $$
B^{-1}=\frac{3}{2\alpha(\eta-\gamma-c_{0})}\left[
\begin{array}{ll}
	\eta&\alpha\\
\frac{\eta+2\gamma+2c_{0}}{3}	& \alpha
\end{array}\right].
$$
Hence $B^{-1}$ is an $M$-matrix, provided that inequality \eqref{Cond} holds. By Theorem \ref{Sta} system \eqref{Ang2} is exponentially stable;
thus system \eqref{Ang} is locally asymptotically stable.\\
\subsection{Stability Analysis for Model 2}

\begin{equation}\label{Ang2tau}
\begin{array}{ll}
\displaystyle\frac{dx}{dt}=\alpha x(t)\ln \left (\frac{K(t)}{x(t)}\right )-p_0x(t) \\[3mm]
\displaystyle\frac{dK}{dt}=(\beta -\gamma ) K(t) -\delta x^{2/3}(h(t)) K(t)- c_0 K(t).
\end{array}
\end{equation}
The system
\begin{equation}
\begin{array}{ll}
\alpha x\ln \left (\frac{K}{x}\right )-p_{0}x=0 \\[3mm]
(\beta-\gamma ) K-\delta x^{2/3} K-c_{0} K=0 \nonumber,
\end{array}
\end{equation}
has a unique positive equilibrium point
\begin{equation} \label{Equi1}
x^{*}=(\frac{\beta-\gamma-c_{0}}{\delta})^{\frac{3}{2}}, ~ K^{*}=x^{*} e^{\frac{p_0}{\alpha}},
\end{equation}
provided that $\beta >\gamma +c_{0} $.\\
Set  $u=\ln\frac{x}{x^{*}}$ and $v=\ln\frac{K}{K^{*}}$, then system \eqref{Ang2tau} has the following form with the trivial equilibrium
 \begin{equation}\label{Ang20}
\begin{array}{ll}
\displaystyle\frac{du}{dt}=-\alpha  u(t)+\alpha v(t)  \\[3mm]
\displaystyle\frac{dv}{dt}= (\beta-\gamma-c_0)(1- e^{^{\frac{2u(h(t))}{3}}}).
\end{array}
\end{equation}

\noindent
Nonlinear system \eqref{Ang20} is equivalent to a second order delay differential equation
\begin{equation}\label{New}
\frac{d^{2}u}{dt^{2}}+\alpha \frac{du}{dt}+A(e^{2u(h(t))/3}-1)=0;
\end{equation}
or, if  $h(t)\equiv t$ , a second order nonlinear ordinary differential equation
\begin{equation}\label{New1}
\frac{d^{2}u}{dt^{2}}+\alpha \frac{du}{dt}+A(e^{2u(t)/3}-1)=0,
\end{equation}
where $ A=\alpha(\beta-\gamma-c_0)$.
Equations \eqref{New}--\eqref{New1} belong to a well-known  class of a Lienard-type  differential equations \cite{Bo1}--\cite{Bu3}.
\begin{equation}\label{ode1}
\frac{d^{2}u}{dt^{2}}+f(u)\frac{du}{dt}+g(u)=0,
\end{equation}
or
\begin{equation}\label{dde1}
\frac{d^{2}u}{dt^{2}}+f(u)\frac{du}{dt}+g(u(h(t)))=0.
\end{equation}
 The following result was obtained by T. Burton \cite{Bu3}.
 \begin{Th}\label{B}
 Suppose functions $f(u)$ and $g(u)$ are continuous with $f(u)>0$ and $ug(u)>0$ if $u\neq 0$.
 Then the zero solution of equation \eqref{ode1} is globally asymptotically stable if and only if
 \begin{equation}  \int_{0}^{\pm \infty}[f(x)+|g(x)|]dx=\pm \infty .\nonumber\end{equation}
 \end{Th}
For equation \eqref{New1} we have
$f(u)=\alpha >0$, $ g(u)=A(e^{2u(t)/3}-1)$, provided $A>0$ and $u(t)>0$ for all $t\geq 0$; thus based on Theorem \ref{B}, system \eqref{New1} is globally asymptotically stable.\\
To prove local stability for Model 2 with delays, we will use the following result recently obtained by L. Berezansky et al. (\cite{Ber2} Corollary 5.2).\\
To this end, consider the linear equation
\begin{equation}\label{lin2}
\frac{d^{2}u}{dt^{2}}+a(t)\frac{du}{dt}+b(t)u(h(t))=0.
\end{equation}
\begin{uess}\label{4}
Suppose $a(t)\geq a_0>0$, $t-h(t)\leq \tau$, and $$\limsup_{t\rightarrow\infty}\frac{|b(t)|}{a(t)}<\tau.$$
Then equation (\ref{lin2}) is exponentially stable.
\end{uess}
\begin{Th}\label{Th2}
Suppose \begin{equation}\label{Op}
t-h(t)\leq \tau ~ and ~ 0<\beta-\gamma-c_0<1.5 \tau.
\end{equation}
Then the positive equilibrium of system (\ref{Ang2tau}) is locally exponentially stable.
\end{Th}
{\bf Proof.}
Stability of the positive equilibrium of system (\ref{Ang2tau}) is equivalent to stability of the trivial solution of equation (\ref{New}). Linearized equation for (\ref{New}) has the form (\ref{lin2}), where $a(t)=\alpha~ and~ b(t)=\frac{2}{3}\alpha(\beta-\gamma-c_{0})$.
Finally, application of Lemma \ref{4} proves the theorem.\\
\begin{Rem} The conditions obtained in Theorem \ref{Th4.1} are independent of delays, whereas local stability of Model 2 (Theorem \ref{Th2}) depends on the  combination of the magnitude of the delays and some parameters of the model. \end{Rem}

\section{ Models with Time-Varying Treatment}
Consider general models (\ref{Ang0}) and (\ref{Ang00}), and assume the existence of the limits
\begin{equation}
\label{limit}
\lim_{t\rightarrow\infty}p(t)=p_0 ~ and ~ \lim_{t\rightarrow\infty}c(t)=c_0.
\end{equation}
We quote the following lemma \cite{HL}.
\begin{uess}\label{L5}
Consider the vector equation
\begin{equation}\label{v1}
\dot{X}(t)=\sum_{k=1}^m A_k(t)X(h_k(t))+F(t),
\end{equation}
where $A_k(t)~ and ~ F(t)$ are  locally essentially bounded matrix and vector functions,
$t-h_k(t)\leq \tau$.
If homogeneous equation
$$
\dot{X}(t)=\sum_{k=1}^m A_k(t)X(h_k(t))
$$
 is exponentially stable, and $\lim_{t\rightarrow\infty}\|F(t)\|=0$,
then for any solution $X$ of equation \eqref{v1}
$$\lim_{t\rightarrow\infty} X(t)=0.$$
\end{uess}
Lemma \ref{L5} is a simple corollary of
the variation of constant formula for solutions of
linear delay differential equations and an exponential
estimation for the fundamental matrix of this equation.\\

Firstly, consider system (\ref{Ang0}) in the form
\begin{equation}\label{Ang_1}
\begin{array}{ll}
\displaystyle\frac{dx}{dt}=\alpha x(t)\ln \left (\frac{K(t)}{x(t)}\right )-p_0x(t)+(p_0-p(t))x(t) \\[3mm]
\displaystyle\frac{dK}{dt}=\beta x(h(t))-\gamma K(t) -\delta x^{2/3}(h(t)) K(t)- c_0 K(t)+(c_0-c(t)) K(t).
\end{array}
\end{equation}

After substitution  $u=\ln\frac{x}{x^{*}}$ and $v=\ln\frac{K}{K^{*}}$ , system (\ref{Ang_1}) has a form
\begin{equation}\label{Ang_2}
\begin{array}{ll}
\displaystyle\frac{du}{dt}= -\alpha  u(t)+\alpha v(t)+p_0-p(t)\\[3mm]
\displaystyle\frac{dv}{dt}=\eta e^{u(h(t))-v(t)} -\gamma -c_0-(\eta-\gamma-c_0 ) e^{^{\frac{2u(h(t))}{3}}}+c_0-c(t),
\end{array}
\end{equation}
where $x^{*} ~and~K^{*}$   are defined by (\ref{Equi}).
Linearization for system (\ref{Ang_2}) yields
\begin{equation}\label{Ang_3}
\begin{array}{ll}
\displaystyle\frac{du}{dt}=-\alpha  u(t)+\alpha v(t) +p_0-p(t) \\[3mm]
\displaystyle\frac{dv}{dt}=\frac{\eta+2\gamma+2c_{0}}{3}u(h(t)) -\eta v(t)+c_0-c(t).
\end{array}
\end{equation}
Lemma \ref{L5} implies the following result.
\begin{Th}\label{Noa1}
Suppose conditions of Theorem \ref{Th4.1} and condition (\ref{limit}) hold.
Then   the pair of numbers $(x^{*}, K^{*})$,  defined by (\ref{Equi}),
is a local attractor for the solutions of system  (\ref{Ang0}):
$$
\lim_{t\rightarrow\infty}x(t)=x^{*}, \lim_{t\rightarrow\infty}K(t)=K^{*}.
$$
\end{Th}
Similar calculations yield
\begin{Th}\label{Noa2}
Suppose conditions \eqref{Op} and (\ref{limit}) hold.
Then   the pair of numbers $(x^{*}, K^{*})$,  defined by (\ref{Equi1}),
is a local attractor for the solutions of system  (\ref{Ang00}).
\end{Th}

\begin{Rem}
Note that the methods applicable for Model 1 yield similar results for the logistic-type model \eqref{AngLog}.
\end{Rem}

\section{Discussion and Open Problems}
Mathematical modeling and simulation can potentially provide insight into the underlying causes of tumor invasion and metastasis, help understand clinical observations, and be of use in designing targeted experiments and assessing treatment strategies \cite{And}, \cite{A} and \cite{EF}.
Compare to some models in literature, where delays were introduced without biological motivation, we give a solid justification for the introduction of the delays into the models.  All  models are perturbed either by a constant therapy or by  time-varying treatments. The existence of unique global solutions for models \eqref{Ang0}--\eqref{AngLog} follows from  Theorem \ref{1}.
For some general class of nonlinear nonautonomous systems with delays we obtained conditions for the existence of a global attractor (Theorem \ref{Sta}). Based on Theorem \ref{Sta}, we proved local stability for Model l. For local stability analysis of Model 2 we used recent result for a Lienard-type second-order differential equation with delays.
Criteria obtained for local attractivity, are explicit and hence are convenient to apply/verify in practice. Note that in Theorems \ref{Noa1}-- \ref{Noa2} the study of nonautonomous equations is reduced to the study of  "limit equations" or asymptotically autonomous systems.
 Note that the analysis of Model 3  could be achieved by using similar techniques.\\
Finally, we formulate some open problems.

\begin{enumerate}
\item
{\bf Conjecture} { Model 1 has a global attractor, provided that a positive equilibrium exists}.
\item {\bf Conjecture} { Model 2 has a global attractor, provided that conditions of the Theorem \ref{B} are satisfied}.
\item {Some additional problems which were not included in the paper: lower and upper estimations of the solutions, extinction, existence and asymptotic stability of periodic solutions,
oscillation and nonoscillation of the solutions about its positive equilibria}.
\end{enumerate}

{\bf   Funding}\\
Research supported in part by: the Israeli Ministry of Absorption and a grant from Vancouver Island University.\\

{\bf Acknowledgements}\\
We wish to express thanks to Dr. A. Gibson (Biology department at Vancouver Island University) whose comments helped us improve the clarity and quality of this paper.

\end{document}